\documentclass[12pt,a4paper,leqno]{article}

\usepackage{fullpage,amsmath,amsthm,amsfonts,amscd,graphicx,
latexsym,amssymb,eucal,curves}
\usepackage[all]{xy}
\newcommand{\ZZ}{\mathbb{Z}}

\newcommand{\CC}{\mathbb{C}}
\newcommand{\PP}{\mathbb{P}}

\newcommand{\HH}{\mathbb{H}}

\newcommand{\Ker}{{\rm Ker}}

\newcommand{\PicS}{{\rm Pic}}
\newcommand{\JacS}{{\rm Jac}}

\newcommand{\Def}{{\rm Def}}

\newcommand{\AlbS}{{\rm Alb}}

\newcommand{\mfS}{{\mathfrak S}}

\newcommand{\mfAV}{{\cal A}}
\newcommand{\mrAV}{A}
\newcommand{\mfC}{{\mathfrak C}}
\newcommand{\mrC}{M}

\newcommand{\Spc}{{\rm Spec \,}}
\newcommand{\echar}{{\chi({\cal O}}}
\newcommand{\sms}{\smallsetminus}

\newcommand{\mi}{{\rm m.i.}}

\def\co{\colon}

\DeclareMathOperator{\pr}{pr}

\newcommand{\ol}{\overline}

\newtheorem{Defi}{Definition}[section]

\newtheorem{Rem}[Defi]{Remark}

\newtheorem{Prop}[Defi]{Proposition}
\newtheorem{Lemma}[Defi]{Lemma}
\newtheorem{Cor}[Defi]{Corollary}
\newtheorem{Thm}[Defi]{Theorem}
\newtheorem{Crit}[Defi]{Criterion}

\begin{document}{\large}
\title{
Maximally irregularly fibred
surfaces of general type}
\author{Martin M{\"o}ller}
\date{}
\maketitle

\begin{quote}
{\footnotesize {\bf Abstract}.
We generalise a method of Xiao Gang to construct 'prototypes'
of fibred surfaces with maximal irregularity without being a product. 
This enables us, in the case of 
fibre genus $g=3$ to describe the possible singular fibres and 
to calculate the invariants of these surfaces. We also prove
structure theorems on the moduli space for fibred surfaces with
fibre genus $g=2$ and $g=3$.

{\bf Mathematics Subject Classification 2000.} 14J10, 14J29, 14D06

{\bf Key words}. fibred surface, moduli space of surfaces of general type, 
high irregularity, fixed part, degenerate fibres}
\end{quote}

\section*{Introduction}
Complex fibred surfaces $f:X \to B$ of small fibre genus 
may be studied by different
techniques according to their irregularity $q(X) = H^1(X,{\cal O}_X)$.
If $b$ denotes the genus of the base $B$ and $g$ the genus of a
fibre, the irregularity of a fibred surface is subject to
$b  \leq q \leq g+b$. In case $q= g+b$ the fibration is trivial, i.e.\
X is birational to a product. 
\newline
Xiao (\cite{Xi85}) and Seiler (\cite{Sei95}) examined surfaces $X$ 
with $g=2$ and irregularity $q(X)=b$ using the fact that these surfaces 
are double coverings of ruled surfaces. Xiao also studies (\cite{Xi85})
surfaces with $g=2$ and $q=b+1$ using the fixed part of the Jacobian
fibration. We extend this technique to surfaces with $q = b+g-1$, 
which we call {\em maximally irregularly fibred}. 
Maximal irregularity implies $g \leq 7$ by (\cite{Xi87})
and maximally irregular fibrations with $g \leq 4$ are 
known to exist (\cite{Pi89}).
\newline
We construct a {\em 'prototype'} for maximally irregularly fibred
surfaces with $g=3$, i.e.\ a fibred surface such that any other 
surface with the same invariants arises via pullback by covering
of the base curves (see Def.~\ref{Prodefi} for the precise definition).
Compared to Xiao's case additional difficulties
arise at the hyperelliptic locus due to the failure of infinitesimal
Torelli. The prototype
enables us to determine the degenerate fibres and invariants
of maximally irregularly fibred surfaces with $g=3$.
\newline
The techniques apply in principle also for
$g=4$ (and, if surfaces exist, also for $g \geq 5$), but these cases 
additionally need an answer to a Schottky type problem, as explained 
at the end of the paper.
\newline
Finally we show that these constructions glue together in
families. We thus obtain structure results for components
of the Gieseker moduli space of surfaces of general type.
The surfaces admitting a maximally irregular fibration form
connected components of the moduli space. These components fibre
over a moduli space of abelian varieties. The fibres are
moduli spaces of stable mappings.
\par
The paper is organised as follows: In \S 1 
we recall some facts on the fixed part of the Jacobian fibration.
We are then able to give the precise definition of prototype.
Omitting some technical conditions on the abelian variety 
$A$ one of the main theorems can then be roughly stated as follows:
\par
{\bf Theorem \ref{goodpro3}' } {\it
Suppose $d \geq 3$ and let $A$ be a $(1,d)$-polarised
abelian variety of dimension $2$.
Then there is a fibred surface $S(A,d) \to B(A,d)$, such
that any maximally irregular fibration $X \to B$ 
with fixed part $A$ and fibre genus $3$ is obtained
via base change $B \to B(A,d)$ from $S(A,d)$.
\newline
The base curve $B(A,d)$ is a double covering of the modular curve $X(d)$.
}
\par
Its proof relies on a parametrisation of the Jacobians of curves into which
the abelian variety $A$ injects. This result is stated as Theorem
\ref{UnivFam} and its proof is presented in \S 2.
\S 3 contains the proof of Theorem \ref{goodpro3} together with the
computation of the invariants of maximally irregular fibrations. 
Finally in \S 4 we use these results to derive some structure results for the
corresponding components of the moduli space.
\newline
Most of the results can be found in the author's thesis (\cite{Mo02}).
\par

\subsubsection*{Acknowledgements}
The author thanks his thesis advisor F.~Herrlich for many
stimulating discussions and a lot of patience. He also thanks
E.~Viehweg for worthful remarks concerning Torelli's theorem.
Some results in the same direction were obtained independently by
J.-X.~Cai. The author thanks him and the referee for his suggestions. 

\subsubsection*{Notation}
{\em We work over the complex numbers throughout}. For a fibred surface
$f: X \to B$ we denote $b= g(B)$ the base genus 
and $g=g(F)$ the fibre genus. The irregularity $q(X) = H^1(X,{\cal O}_X)$
is subject to $b \leq q \leq b+g$. We say the fibration
$f$ is {\em of type $(g,b)$}.
\newline
Let $\mfC_g(\cdot)$ (resp.\ $\mfAV_g(\cdot)$) denote the moduli functor for
smooth curves of genus $g$ (resp.\ abelian varieties of dimension $g$)
and let $\mrC_g$ (resp.\ $\mrAV_g)$) denote the corresponding
coarse moduli spaces.
\newline 
When discussing surfaces of general type, we denote by $X$ its {\em canonical
model}, i.e.\ a normal surface with $K_X$ ample and at most
rational double points. If necessary, $S$ denotes the corresponding
smooth minimal model. A {\em relative canonical model} of a flat
family of surfaces of general type over some base $T$ has 
a canonical model in the fibre over each complex point of $T$. 
Each flat family of surfaces of general type is birational to
a relative canonical model (see \cite{Tv72}). 
We denote by $\mfS(\cdot)$ 
the moduli functor which associates with a scheme $T$
the set of flat families of relative canonical models of surfaces 
of general type over $T$.   

\section{Prototypes for fibred surfaces} \label{Protosec}

A {\em fibration} of a surface $f:X \to B$ is a surjection onto a smooth
curve $B$ with connected fibres. The fibration is called
{\em regular} if $q(X) = b$ and {\em irregular otherwise}.
If $q(X) = b + g$, the surface $X$ is birational to a product
of the base curve and smooth fibre. 
If the moduli map $B \to M_g$ has image reduced to a
point the fibration is said to be {\em isotrivial} or
{\em of constant moduli}. 
\newline
In the sequel we are interested in the case of non-trivial fibrations
with maximal irre\-gularity, i.e.\ with $q(X) = g+b-1$. We call them
{\em maximally irregularly fibred}. The Jacobians of the fibres of
such a surface have a large abelian variety in common, the
{\em fixed part} ('partie fixe' or '$L/K$-trace' in \cite{La59}).
We recall some facts about the fixed part of a  fibred family 
$X \to B \to T$ of surfaces. The reader less interested in 
statements on the moduli space may think of $T = \Spc \CC$ 
in sections $1-3$.
\par
Let $A = \PicS^0_{X/T}/\JacS_{B/T}$.
Over the locus $B'$ where $X \to B$ is smooth (which is dense
in every fibre over $T$) $\JacS_{X'/B'}$ exists. The Picard 
functor applied to the $B'$-morphism $X' \to X \times_T B'$ induces
$$ i'_A: \PicS^0_{X \times_T B'/B'} \to \JacS_{X'/B'},$$
which factors through an injection
$$i_A: A\times_T B' \hookrightarrow \JacS_{X'/B'}. $$
\par
We want to provide $A$ with a polarisation. To that
purpose we restrict the principal polarisation $\Theta_J$ of
$\JacS_{X'/B'}$ to $A \times_T B'$. We consider polarisations
as (relatively) ample line bundles, but keep the classical
divisor notation. By the rigidity theorem
the polarisation comes from a relatively ample line bundle
$\Theta_A$ on $A$, i.e.\ $i_A^* \Theta_J = \pr_A^* (\Theta_A)$,
where $\pr_A: A \times_T B' \to A$ is the projection.
\par
\begin{Defi} \label{assdegdeg}
The polarised abelian scheme $(A/T, \Theta_A)$ is called
the {\em fixed part} of the fibration $X \to B \to T$. If the
degree of isogeny $\lambda(\Theta_A)$ is $d^2$, we call $d$ 
the {\em associated degree} of the fibred family. We also
say that $X \to B \to T$ is {\em of type $(A/T, \Theta_A)$}.
\end{Defi}
\par
\begin{Rem} \label{AlbRem}
{\rm i) In \cite{Xi92}, Xiao calls the degree of
$\alpha: X \to \AlbS_X$ (if generically finite) 'associated
degree'. We prefer to call $\deg(\alpha)$ the {\em Albanese degree} (denoted
by $\gamma$). 
\newline
ii) We can obviously extend the notion 
of associated degree to any injection 
$ A \to J$ of abelian varieties, whenever $J$ carries a principal
polarisation. 
}\end{Rem}
\par
\begin{Defi} \label{Prodefi} Let $(A,\Theta_A)$ be a polarised
abelian variety. 
A relatively minimal model of a fibred surface
$S(A,d) \to B(A,d)$ is called a
{\em prototype for fibred surfaces with fixed part  $(A, \Theta_A)$,
if each fibred family of surfaces $X \to B \to T \in \mfS(T)$ of
type $(g,b)$ with associated degree $d$ and fixed part 
$(A \times T, \pr_1^* \Theta_A)$}
is the relative canonical model of
$B \times_{B(A,d)} S(A,d) \to B$
for a surjection $\varphi: B \to B(A,d)$. The latter is called a {\em
prototype base change}.
\par
If conversely the pullback by any such $\varphi$ gives a fibred
surface with fixed part $A$, we call $S(A,d) \to B(A,d)$ a {\em 
good prototype}.
\end{Defi}
\par
We do not demand $\varphi$ to be unique. However $\varphi$
turns out to be more or less unique in the cases studied below.
\par

\subsection{Prototypes for maximally irregular fibrations
with $d \geq 3$ and $g=2$ or $g=3$}

Let $X(d)$ be the completion of the modular curve $X'(d) = \HH/\Gamma(d)$, 
where $\HH$ is the upper half plane and $\Gamma(d)$ the principal
congruence subgroup of level $d$. Then we have:
\par
\begin{Thm} \label{goodpro3}
Let $(A,\Theta_A)$ be a $(1,d)$-polarised
abelian variety of dimension $2$ without a nontrivial principally 
polarised abelian subvariety and suppose $d \geq 3$.
\newline
Then there is a fibred surface $S(A,d) \to B(A,d)$, which
is a good prototype for maximally irregular fibrations
with fixed part $(A, \Theta_A)$ and fibre genus $3$.
The base curve $B(A,d)$ is a double covering of the modular curve $X(d)$.
\newline
Moreover, the prototype base change $\varphi$ is unique up to composition
with the involution $\sigma$ of $B(A,d)$ over $X(d)$.
\end{Thm}
\par
The restrictions on $A$ are necessary: If $A$ had a nontrivial 
principally polarised abelian subvariety, the
Jacobian of each fibre would be reducible, a contradiction. And
in Section \ref{secAbVar} we will see that the polarisation of the fixed
part in case of maximal irregularity is always of type
$(1,\ldots,1,d)$. Here and everywhere in the sequel $A$ should be
considered as a polarised variety. In particular the fibred surface 
$S(A,d)$ depends on the polarisation. To avoid overloading notations we
frequently omit the polarisation.
\par
This theorem should be compared with the following theorem, which is
basically due to Xiao.
\par
\begin{Thm} \label{goodpro2}
For each one-dimensional abelian variety $A$ and each $d \geq 3$
there is a fibred surface $S(A,d) \to B(A,d)$, which is a good
prototype for maximally irregular fibrations of type $(2,b)$ with
fixed part $A$ and associated degree $d$. The base curve $B(A,d)$
is the modular curve $X(d)$, it does not depend on $A$.
\newline
Moreover the prototype base change $\varphi$ is unique.
\end{Thm}
\par
The proofs of these theorems will rely on the following parametrisation
of injections of abelian varieties. The reader should think of 
the injection of the fixed part of the relative Jacobian of a
fibred surface into the Jacobian of a fibre. $\pr_A: A \times X \to A$ always 
denotes the first projection.
\par
\begin{Thm} \label{UnivFam}
Let $(A,\Theta_A)$ be a $(g\!-\!1)$ dimensional abelian variety
with a polarisation of type $\delta = (1,\ldots,1,d)$. 
If $d \geq 3$ there is a family 
$$ (j'(A,d): J'(A,d) \to X'(d), \quad \Theta_{J'(A,d)})$$
of $g$-dimensional, principally polarised
abelian varieties with an injection 
$$i_{X'(d)}: A \times X'(d) \to J'(A,d)$$
over $X'(d)$,
such that $i_{X'(d)}^* \Theta_{J'(A,d)} = pr_A^* \Theta_A$. The 
family $j'(A,d)$ is universal
for principally polarised abelian varieties $(J/T, \Theta_J)$ with
an injection $i_T: A \times T \to J$ over $T$, such that $i_T^*(\Theta_J) = 
\pr_A^* \Theta_A$.
\end{Thm}
\par
In section \ref{ProtandInv} we will obtain the prototypes by applying a
Torelli theorem to this universal family. As the Torelli map has no
longer dense image in $\mrAV_g$ for $g>3$, our strategy is limited
to $g \leq 3$. A relative version
of Thm.\ \ref{UnivFam} also holds true and will be needed in section
\ref{secMSpace}.
\par

\subsection{Maximally irregular fibrations
with $d=2$ and $g =2$ or $g=3$, Isotriviality}

The case $d=2$ needs a separated treatment because 
$(-1)$ does not act freely on $\HH$ and hence Thm.\ \ref{UnivFam}
is false for $d=2$.
For $g=2$ and $d=2$ Xiao showed (\cite{Xi85} Example 3.1) 
that these surfaces are double coverings of principal 
homogeneous spaces for $E \times B$, where $E$ is the 
fixed part i.e. an elliptic curve. 
\newline
We analyse now the case $g=3$ and the possibility of
a maximally irregular fibration to be isotrivial or
to have a hyperelliptic generic fibre.
\par
\begin{Prop} \label{genHELok}
Let $f: X\to B$ be a maximally irregular fibration of type
$(g,b)$ with $g \geq 3$. If the generic fibre is
hyperelliptic or if $X \to B$ is isotrivial, 
we have $g=3$, $d=2$ and the Albanese degree $\gamma$  
equals $2$.
\newline 
Conversely if $f: X \to B$ has $g=3$ and $\gamma=2$, then $f$
is isotrivial.
\end{Prop}
\par
{\bf Proof:}
From \cite{Pi89} or \cite{Xi92} it follows that the image
of the Albanese map of a hyperelliptic fibration is a product
of the base curve and a curve $C$ of genus $g-1$. Riemann-Hurwitz
implies that $g=3$ and that each fibre is an unramified double 
cover of $C$, hence $\gamma=2$. By \cite{LB92} Theorem 12.3.3 
we have $d=2$.
If the fibration is isotrivial, Prop.\ 2.2 of \cite{Ser96} 
implies that the fibre of the Albanese image has dimension $g-1$
and we conclude as above.
\newline
For the converse note that $\gamma =2$ implies that each fibre
of $f$ is an unramified double cover of the same curve of
genus $2$. 
\hfill $\Box$
\par
See \cite{Ser96} or \cite{Ca00} for more on isotrivial
fibrations, which we exclude from now on. For section \ref{secMSpace}
we note that isotrivial fibrations of maximal irregularity 
form components of the moduli space,
because $\gamma$ is constant on connected components.
\newline
We will fix the 
$4$ complements $\Gamma(2)_S$ of $\pm 1$ in $\Gamma(2)$. Thereby the index
$S$ denotes the set of irregular cusps (see \cite{Sh71})
of $\HH/\Gamma(2)_S$, i.e.\ $S$ is a subset of $\{0,1, \infty\}$
of order $1$ or $3$.
\newline
Let $(A,\Theta_A)$ be a $(1,2)$-polarised abelian
surface without non-trivial principal polarised abelian
subvariety. We thus obtain a (certainly not good) prototype also in
this case:
\par
\begin{Thm} \label{3protod2}
For each $S$ as above there is a fibred surface
$$h(A,2)_S: S(A,2)_S \to \PP^1$$
which is a prototype in the following sense:
\newline
For each family of surfaces $X \to B \to T$ of type $(3,b)$
with maximal irregularity, $d=2$, $\gamma =1$ and of type 
$(A, \Theta_A)$, there is an index set $S$ and a unique 
morphism $\varphi: B \to \PP^1$, such that $X$ is the canonical model of
$S(A,2)_S \times_{\PP^1} B$.
\end{Thm}
\par

\section{Moduli spaces for abelian varieties with a fixed
injection}
\label{secAbVar}

We fix an injection of abelian varieties $i_A:A \to J$,
where $J$ is principally polarised and we consider the
case $\dim J = g$, and $\dim A = g\!-\!1$. Then there
is a complementary abelian subvariety 
$E=\Ker(i^{\vee}:J \to A^{\vee}$) of dimension one, where $A^{\vee}$ is
the dual abelian variety. This kernel is indeed connected
(see \cite{LB92} Section 12.1) and we denote by $i_E: E \to J$ 
the inclusion. By loc.\ cit.\
Cor.\ 12.1.5 there is a $d$ such that the
polarisations $\Theta_A=i_A^* \Theta_J $ and $\Theta_E = i_E^* \Theta_J$
are of type $(1,\ldots,1,d)$ and $(d)$ respectively. If $i_A$
comes from a fibred surface, $d$ coincides with the associated
degree in Def.\ \ref{assdegdeg}.
\newline
Let $V = H^0(J, \Omega_J)^{\vee}$ and $U = H_1(J,\ZZ)$ be
the uniformisation of $J = V/U$. Further let $V_A = H^0(A, \Omega_A)^{\vee}$
and $V_E=H^0(E, \Omega_E)^{\vee}$ and consider them as subspaces
of $V$ via $i_A$ and $i_E$. We then obtain
$U_A := H_1(A, \ZZ) = V_A \cap U$ and $U_E := H_1(E,\ZZ) = V_E \cap U$.
\par
\begin{Lemma} \label{Hombaprop}
In the above situation, there exists a
basis of homology adapted to the injections $i_A$ and
$i_E$ ('adapted basis' for short)
$$B = \{ u_1, \ldots, u_{g}, u_{g+1}, \ldots, u_{2g-2}, u_{2g+1}, u_{2g+2} 
\}$$ of $U$, 
i.e. a basis, such that        
$$ u_1, \ldots u_{g-1}, u_{g+1}, \ldots,  u_{2g-1} = d u_{2g+1} - u_{g}$$ 
is a symplectic basis of $U_A$ and
$$ u_{g}, u_{2g} = d u_{2g+2} - u_{g-1} $$
is a symplectic basis of $U_E$.
\end{Lemma}
\par
\par
{\bf Proof:} We denote the alternating form on $U$ induced by
$\Theta_J$ by $(\cdot,\cdot)$. The polarisations $\Theta_A$
and $\Theta_E$ induce alternating forms, which are by construction
the restriction of $(\cdot,\cdot)$ to $U_A$ and $U_E$.
Denote the elements of a symplectic basis of $U_A \oplus
U_E$ of type $(1,\ldots,1,d,d)$ by $u_1,\ldots,u_{g-2},a_1,a_2,
u_{g+1},\ldots,u_{2g-2},b_1,b_2$. The projection $p: V \to V/V_A$
induces an isomorphism $U/(U_A\oplus U_E) \to p(U)/p(U_E) \cong
\ZZ/d_1\ZZ \oplus \ZZ/ d_2\ZZ$. Let $v_1, v_2 \in U$ be generators
of this quotient. The type of the restricted alternating form
implies that $d$ divides all of $(d_1 v_1,a_i)$, $(d_1 v_1, b_i)$,
$(d_2 v_2, a_i)$ and $(d_2 v_2, b_i)$ for $i=1,2$. As neither
$v_1$ nor $v_2$ is orthogonal to $\langle a_1,a_2,b_1,b_2\rangle$
this is only possible if $d=d_1=d_2$.
\newline
Hence we find $a,a' \in U_A$ and $b, b' \in U_E$ such that
$v_1 = (a + b)/d$ and $v_2 = (a' + b')/d$.
The images $p(v_1), p(v_2)$ generate $p(U)/p(U_E)$ and
thus $\{b, b' \}$ is, changing the order if necessary,
a symplectic basis of $U_E$. Now $U_A$ and $U_E$ are orthogonal
with respect to $(\cdot,\cdot)$ and 
$\{u_1,\ldots,u_{g-2}, u_{g+1},\ldots,u_{2g-2},a,a',b,b'\}$
is also a basis of $U_A \oplus U_E$. The determinant of $(\cdot,\cdot)$
restricted to $U_A \oplus U_E$ is $d^2$ and together this implies
$(a,a') = d$. Hence we can take
$u_{g-1} = a'$, $u_g = b$, $u_{2g+1}=v_1$ and $u_{2g+2}=v_2$.
\hfill $\Box$
\par 
A similarly constructed basis for $g=2$ is 
called 'base homologu\'ee' in \cite{Xi85}.
\par
We now check how far this basis of homology is from being unique.
\par
\begin{Lemma} \label{HomobasisW}
Let $B$ be an adapted basis of $U$ and let
$$ (u_{g}', u_{2g}')^T = M (u_{g}, u_{2g})^T, \quad \quad
M \in SL_2(\ZZ).$$
If we let $u_{2g+1}' =\frac{1}{d}(u_g' + u_{2g-1}) $ and 
$u_{2g+2}' =\frac{1}{d}(u_{2g}' +u_{g-1}) $, then
$$ B' = \{ u_1, \ldots, u_{g-1}, u_{g}', u_{g+1}, \ldots, u_{2g-2}, 
u_{2g+1}', u_{2g+2}' \}$$ 
is another adapted basis $U$, if and only if 
$$ M \in \Gamma(d) =  
 \{ M \in SL_2(\ZZ), \quad M \equiv \left(\begin{array}{rr}
1 & 0 \\
0 & 1 \\
\end{array} \right) \! \mod d \}.$$
\end{Lemma}
\par
{\bf Proof:} This change of basis is independent of the
elements $u_1, \ldots, u_{g-2}, u_{g+1}, \ldots
u_{2g-2}$. We can hence copy the proof of
Lemme 3.7 in \cite{Xi85} literally.
\hfill $\Box$
\par
We now start the construction the family $j'(A,d)$ 
of Thm.\ \ref{UnivFam}.
With respect to a symplectic basis 
the period matrix of $A$ is of the form
 $$\left(
\begin{matrix}
 z_{11} & \ldots & z_{1\, g-1} & 1 &\ldots& 0 & 0 \\
\vdots && \vdots && \ddots \\
 z_{g-1\, 1} & \ldots & z_{g-1\, g-1} & 0 & \ldots &0 &d \\
\end{matrix}
\right)^t  = (Z, \rm{diag}(\delta))^t,$$
where $\delta= (1,\ldots,1,d)$, for some $Z$ in the
Siegel half space $\HH_{g-1}$. 
Let $$p\co\  {\cal V}\co=\HH \times \CC^g\to \HH$$
be the trivial vector bundle over the upper half plane $\HH$. For $z \in \HH$
the sections 
$$\begin{array}{lcllcll}
u_r(z) &=& \phantom{\frac{1}{d}}(z_{r1},\ldots,z_{r\,g-1}, 0)& 
u_{g+r}(z) &=& \phantom{\frac{1}{d}}(e_r, 0) & 
\hbox{for}\, r=1,\ldots,g-1 \\
u_{g}(z) &=& \phantom{\frac{1}{d}}(0,\ldots,0,0,z) & 
u_{2g}(z) &=&\phantom{\frac{1}{d}}(0_{g-1},1) \\
u_{2g+1}(z)&=&\frac{1}{d}(0,\ldots,0,d,z) &
u_{2g+2}(z) &=&\frac{1}{d}(z_{g-1}, d)\\
\end{array} $$
define a lattice $U(z)$ in ${\cal V}$,
where $e_r$ are the rows of ${\rm diag}(\delta)$.
The quotients $J(z) = {\cal V}/U(z)$ are a family $J(z)\to \HH$ 
of complex tori. They admit an  alternating form
of type $\delta' = (1,\ldots,1,d,d)$ 
with respect to ${u_1, \ldots, u_{2g}}$.
This defines a principal polarisation on $J(z)$, which
is thus a family of abelian varieties. The injection of
$V_A = \CC^{g-1}$ into the first $g-1$ components
and the injection of the lattices 
$$\langle u_1(z),\ldots,u_{g-1}(z),u_{g+1}(z),\ldots,u_{2g-1}(z)
\rangle \hookrightarrow U(z)$$
defines an injection of the trivial family $A \times \HH$
into $J(z)$ and by construction the restriction of the 
principal polarisation to this lattice is of the right type.
Due to Lemma \ref{HomobasisW} we can define a
$\Gamma(d)$-action on this family, which respects
this injection:
$M = \left(\begin{array}{cc}
\alpha & \beta \\
\gamma & \delta \\ \end{array}\right)  \in \Gamma(d)$
acts on $\HH$ as usual and we let
$$ \begin{array}{lcll}
Mu_i & = &\phantom{-} u_i &\hbox{for}\, i=1,\ldots,g\!-\!1,g\!+\!1,
\ldots,2g\!-\!1 \\ 
Mu_g(z) &= &\phantom{-}\delta u_g(Mz) -\beta u_{2g}(Mz) 
& = (0_{g-1}, z/\gamma z + \delta) \\
Mu_{2g}(z)& =& -\gamma u_g(Mz) + \alpha u_{2g}(Mz)
& = (0_{g-1}, 1/\gamma z + \delta). \\
\end{array}$$
By letting $u_{2g+1} = \frac{1}{d}(u_{2g-1} + u_g)$
and $u_{2g+2} = \frac{1}{d}(u_{g-1} + u_{2g})$
we have defined a $\Gamma(d)$-action on $U$  
such that $M: J(z) \to J(Mz)$ is an isomorphism over $A$.
\par
{\bf Proof of Thm.\ \ref{UnivFam}}: For $d \geq 3$
the group $\Gamma(d)$ acts freely on $\HH$, so by taking quotients 
we obtain the desired family $$ j'(A,d): J'(A,d) \to X'(d)$$
of principally polarised abelian varieties.
It remains to show the universal property. 
\newline
Let $(j:J \to T, \Theta_J)$ and $i_A$ be as in the statement of the
theorem and let $\pi: E \to T$ with $i_E: E \to J$ be the
family of complementary abelian varieties.
Let $R_1 \pi_* \ZZ$ denote the local system in the complex
topology with fibres $H_1(E_t,\ZZ)$. Over a sufficiently small 
complex chart $T_i$ the $R_1 \pi_* \ZZ$ and $R_1 j_* \ZZ$ are free. 
Hence we can find sections $u^i_{g}$ and $u^i_{2g}$ in 
$(R_1 \pi_* \ZZ)(T_i)$ that form a symplectic basis which
can be completed in $(R_1 j_* \ZZ)(T_i)$ to sections $u^i_1,\ldots
u^i_{2g+2}$  with the above properties. We can map these
sections to $(j_* \Omega_J)^{\vee}$ and project
them modulo $i_A(({\pr_i}_* \Omega_A)^{\vee})$
(where $\pr_i: A\times T_i \to T_i$ is the projection).
The quotient of the two projected sections defines a map
$\varphi_i: T_i \to \HH$.
Given another chart $T_k$ we can find sections
$u^k_{g}$ and $u^k_{2g}$ which define $\varphi_k:T_k \to \HH$.
By the above lemma the morphisms $T_i \to \HH/\Gamma(d)$
glue to the desired morphism $\varphi: T \to \HH/\Gamma(d)$. 
\newline
We finally remark, that all these complex spaces are in fact
algebraic, see \cite{LB92} Section 10.8.
\hfill $\Box$
\par
We will now calculate the monodromy of $j'(A,d)$. 
For that purpose we fix the following symplectic
basis of $H_1(J,\ZZ)$ (where $J$ is a fibre of $j'(A,d)$):
$$ \alpha_r = u_r \;(r=1,\ldots,g\!-\!2),\; \alpha_{g-1} = u_{2g+2},\; 
\alpha_g = u_{2g+1},\; \beta_r = u_{g+r} \; (r=1,\ldots,g)  $$ 
The period matrix $T(z)$ of $J$ with respect to this basis is
$$ 
T(z) = 
\left(
\begin{array}{cc}
Z  & Z_{12} \\
Z_{21} & z \\
\end{array}
\right), \quad \text{where} \; Z_{21} = (0,\ldots,0,1), \; Z_{12} = Z_{21}^T. 
$$
\par
\begin{Lemma} \label{MonodrLemma}
The monodromy of $j'(A,d)$ along a path $\gamma$ around $\infty$ is
$$ M(\gamma)  = I_{2g} + E_{g,2g} $$
\end{Lemma}
\par
{\bf Proof:} 
The path $\gamma$ lifts to a path from $z$ to $z+d$ in $\HH$. The result
follows from noting that
$M(\gamma$) satisfies $T(z+d) = M(\gamma)\cdot T(z)$ (see \cite{Na74}).
\hfill $\Box$
\par
\begin{Rem} {\rm 
As the monodromy is unipotent, \cite{FC90} Th.\ V.6.7 shows that there
is a unique extension of $j'(A,d)$ 
to a semi-abelian scheme over $X(d)$, which we denote by
$ j(A,d): J(A,d) \to X(d)$.
}\end{Rem}
\par
\begin{Rem} \label{mrem} {\rm 
The choice of a symplectic basis of $H_1(J,\ZZ)$ as above
defines an injection $\phi: \Gamma(d) \to Sp(g,\ZZ)$ into the
symplectic group and an injection $\HH \to \HH_g$ equivariant with respect to
$\phi$ and the natural actions of these groups. The moduli map
of the family $j'(A,d)$ is hence an immersion
$$m': X'(d) \to A_g.$$ 
\newline
To prove the prototype theorems, we need a version of this map
using fine moduli spaces. Fix a level-$[n]$-structure
on the abelian varieties parametrised by $X(d)$.
Writing down $\phi$ explicitely, we
see that a level-$[n]$-structure on $J(z) \to \HH$
is invariant under $M \in \Gamma(d)$, if $M$ is actually 
in $\Gamma(nd)$. Let $X'_n(d)$ denote $\HH/\Gamma(nd)$. This curve equals
of course $X'(nd)$, we just want to emphasize its different role here.
Let $j'_n(A,d): J'_n(A,d) \to X'_n(d)$ be the pullback of $J'(A,d)$.
$X'_n(d)$ is the moduli space of abelian threefolds with an
injection by $A$ plus the level-$[n]$-structure. Now
the moduli map of $j'_n(A,d)$ is an immersion
$$m'_n: X'_n(d) \to \mrAV_g^{[n]} $$
we were heading for. Note that both $m'$ and $m'_n$ depend on $A$. 
}\end{Rem}
\par
\begin{Lemma} \label{3noredlemma}
For $g=3$ and $A$ without non-trivial principally polarised
subvariety almost all fibres of $j'(A,d)$ are not reducible as abelian
varieties with polarisation.
\newline
If moreover $A$ does not contain any elliptic curve then none
of the fibres of $j'(A,d)$ is reducible as abelian
variety with polarisation. 
\end{Lemma}
\par
{\bf Proof:}
We first prove the second statement: Suppose the contrary is the 
case for $t \in X(d)$, i.e.\ 
$$(J_t, \Theta_{J_t}) \cong (A'\times E', p_1^* \Theta_{A'} + p_2
\Theta_{E'}).$$ 
By hypothesis $p_2 \circ i_A: A \to E'$ has to
be the zero map. Hence the injectivity of $i_A$ implies that
it is actually an isomorphism. But then $i_A^* \Theta_A$ 
is a principal polarisation, not of type $(1,d)$.
\par
Suppose under the hypothesis of the first statement that
we have for $t \in X(d)$ a splitting as above.
Consider the inclusion $i_{E_t}: E_t \to J_t$ of the
complementary abelian variety. If $p_2 \circ i_{E_t}: E_t \to E'$
was zero the image of $i_A^\vee|_{E'}$ would be a principally
polarized subvariety of $A$. Otherwise $i_{E_t}^\vee p_2 \circ i_A$
gives an isogeny of an elliptic curve in $A$ to $E_t$. This
cannot happen for all $t\in X(d)$. Hence it happens only 
a finite number of times.
\hfill $\Box$
\par
\begin{Rem} {\rm 
If $A$ contains an elliptic curve $E$ there may be
fibres of $j'(A,d)$, that are reducible with
polarisation.
\newline
To see this take $(A,\Theta_0)$ principally polarised, 
containing an elliptic curve $E$ but irreducible
with polarisation. Let $d' = \deg (\Theta_0|_E)$. Provide
$E$ with a principal polarisation $\Theta_E$ and let
$i^\vee: A \to E$ be the dual of the inclusion. Then
$$ (id_A \times i^\vee): A \to J := A \times E $$
is injective and the pullback of the principal polarisation 
$p_1^* \Theta_0 \otimes p_2^*\Theta_E$ is a polarization of 
type $(1,d'+1)$ on $A$.  
}
\end{Rem}

\section{Proof of the prototype theorems} \label{ProtandInv}

To prove Thm.\ \ref{goodpro3} we want to apply a Torelli theorem to
the family constructed in Thm.~\ref{UnivFam}. The fact that
the map $i:\mrC_g \to \mrAV_g$ from the moduli space of curves
to the moduli space of abelian varieties is an isomorphism (see
\cite{OS80}) is not sufficient here, because we have to pull back 
the universal families. We therefore use auxiliary level-$[n]$-structures.
Let $g \geq 3$, $n\geq 3$ and 
$\Sigma$ be the involution 
on $\mrC_g^{[n]}$, sending the level structure $\alpha$ to $-\alpha$. 
We call $V_g^{[n]}$ the quotient by this involution.
$$\xymatrix{
\mrC_g^{[n]} \ar[d]^q \ar[dr]^{i^{[n]}} &\\
V_g^{[n]} \ar[d] \ar[r]^{i_V^{[n]}} & \mrAV_g^{[n]} \ar[d] \\
\mrC_g \ar[r]^i & \mrAV_g \\
} $$
Let $H \subset \mrC_g$ be the hyperelliptic locus and
$H_V$ its preimage in $V_g^{[n]}$. The map
$q$ is ramified exactly over $H_V$; over $V_g^{[n]} \sms H_V$
there is a universal family of curves. 
\par
\begin{Thm} (\cite{OS80} Thm.\ 3.1) The map
$i_V^{[n]}$ is an embedding.
\end{Thm}
\par
\begin{Prop}
Let $h_i: C_i \to T$ for $i=1,2$ be families of smooth curves
of genus $g \geq 2$ without hyperelliptic fibres, with
sections $s_i$ and induced embeddings $f_i: C_i \to \JacS_{C_i/T}$. 
If $\beta: \JacS_{C_1/T} \to \JacS_{C_2/T}$ is an isomorphism,
there exists a unique isomorphism $\alpha: C_1 \to C_2$ and
a translation $t_c$, such that
$$f_2 \circ \alpha = t_c \circ \beta \circ f_1.$$
\end{Prop}
\par
{\bf Proof:} For $T = \Spc \CC$ this is 
\cite{Mi86} Theorem 12.1.  If $T$ is the spectrum a local Artinian ring, 
\cite{OS80} Proposition 2.5 implies the existence of the isomorphism,
which is unique on the special fibre and therefore unique, because
curves of genus $\geq 2$ do not have infinitesimal automorphisms.
The general statement now follows thanks to the uniqueness by
descent.
\hfill $\Box$
\par
\begin{Cor} \label{JtoC}
Let $j:J \to T$ be a principally polarised 
abelian scheme of dimension $g \geq 3$, such that
the image of the induced mapping $\varphi:T \to \mrAV_g$ 
is contained in $i(\mrC_g \sms H)$. Then there is a unique
family of curves $h: C \to T$, whose Jacobian is, after an
appropriate faithfully flat base change, isomorphic to $J$.
\end{Cor}
\par
{\bf Proof:} We make a base change $\tilde{T} \to T$ such
that $\tilde{j}: \tilde{J} = J \times_T \tilde{T} \to \tilde{T}$ 
admits a level-$[n]$-structure
in order to use the above results.
Fixing a level structure and letting $\tilde{\varphi}:
\tilde{T} \to \mrAV_g^{[n]}$, we obtain a family of curves
$\tilde{h}: \tilde{C} \to \tilde{T}$ by pulling
back the universal family over $V_g^{[n]} \sms H_V$ via
$(i_V^{[n]})^{-1} \circ \tilde{\varphi}$.
Fix an isomorphism $\JacS_{\tilde{C}/\tilde{T}} \to \tilde{J}$ (which is
unique only up to $\pm 1$) and suppose $\tilde{h}$ admits
a section $\tilde{s}$. We can assure this by making another base change.
$\tilde{s}$ induces an embedding $\tilde{f}: \tilde{C} \to \JacS_{\tilde{C}/
\tilde{T}}$. The natural descent data on $\tilde{J}$ give
an isomorphism $\beta_{12}: \pr_1^* \JacS_{\tilde{C}/\tilde{T}} 
\to \pr_2^* \JacS_{\tilde{C}/\tilde{T}}$, where as usual $\pr_i$
denotes the projections $\tilde{T} \times_T \tilde{T} \to \tilde{T}$.
By the above proposition we obtain an isomorphism $\alpha_{12}$ 
between the pullbacks of the families of curves, commuting
with $\beta_{12}$ up to translation.
It satisfies the cocycle condition, because
the descent data on $\tilde{J}$ do so.
As curves of genus $\geq 2$ come along with an ample
canonical sheaf, the descent data are effective (see e.g.
\cite{BLR90} Theorem 6.7).
\newline
It remains to show uniqueness: 
Let $C_i \to T$ for $i=1,2$ be two families of curves
with the desired property. Again after base change
$\tilde{T} \to T$ we have isomorphisms
$\lambda_i: \JacS_{\tilde{C_i}/\tilde{T}} \to \tilde{J}$,
unique up to $\pm 1$. Fixing them, the cocycle condition
on $\JacS_{\tilde{C_i}/\tilde{T}}$ and
$\tilde{J}$ implies that $\lambda_2^{-1} \circ \lambda_1$
respects the descent data. By the above proposition this
gives a unique isomorphism between $\tilde{C_1}$ and $\tilde{C_2}$,
commuting with the embeddings and $\lambda_2^{-1} \circ \lambda_1$ 
up to translation. Hence this isomorphism descends to
the one between $C_1 \to C_2$ we sought.
\hfill $\Box$
\par
\begin{Rem}{\rm 
Letting $p: \tilde{T} \times_T \tilde{T} \to T$ 
the morphism $\beta_{12}$ in the above proof is induced by 
the morphism $p^* C \to p^* C$ coming from the exchange of
factors only up to $\pm 1$ . Lemma \ref{HEdurchnd} 
gives an example where $\beta_{12}$ does indeed not
coincide with this morphism. And this Lemma also
reveals that we cannot replace 'faithfully flat' by '{\'e}tale' in
the above corollary.} \end{Rem}
\par
{\bf Proof of Thm.\ \ref{goodpro3}:} 
We check (Step 1) that the image of the moduli map
$m'$ lies generically in the image of $i(\mrC_g \sms H)$.
We then (Step 2) apply Cor.\ \ref{JtoC} to pullback
the family $j(d)$ via $i$ and complete it to a prototype
${\cal S}(A,d) \to X(d)$. We show that the prototype
base change $\varphi$ is unique.
In Step 3 we provide all families with auxiliary level-$[n]$-structures.
The double covering $\mrC_g^{[n]} \to V_g^{[n]}$ will induces
a double covering of $X'_n(d)$. When we finally divide out
the level-structures we obtain the double covering $B(A,d) \to X(d)$
and a good prototype $S(A,d) \to B(A,d)$. It will turn out
in Cor.\ \ref{covramhyp} that the double covering is ramified, in particular
connected, hence that ${\cal S}(A,d) \to X(d)$ was indeed not a good prototype.
For simplification we omit the dependence on $(A,d)$ from notation. 
\par
{\it Step 1:} The image of the moduli map $m'_n: X'_n(d) \to \mrAV_g^{[n]}$
lies generically in the image of $i_V^{[n]}$ by Lemma \ref{3noredlemma}.
Let $B'_n \to X_n(d)' \times_{V_3^{[n]}} \mrC_3^{[n]}$ be
the normalisation of the fibre product, in which the first map is
$(i_V^{[n]})^{-1} \circ m'_n$ and let $B_n$ be the smooth
completion of $B'_n$. If $m'(X'(d))$ was contained
in the image $i(H)$ of the hyperelliptic locus, the pullback
of the universal family over $\mrC_3^{[n]}$ to $B'_n$
would have hyperelliptic generic fibre. By construction its
completion to a fibred surface over $B_n$ has maximal irregularity
and $d \geq 3$.  This contradicts Prop.\ \ref{genHELok}.
\par 
{\it Step 2:} We now apply Cor.\ \ref{JtoC} and 
obtain a family of curves over the non-hyperelliptic locus
of $X'(d)$. The relatively minimal model of a completion is denoted
by $h: {\cal S} \to X(d)$.
\newline
To prove the universal property of the prototype, let $X \to C \to T$
be a fibred family of surfaces with fixed part $A \times T$. Let
$C'$ the locus with smooth fibres.
As explained in Section \ref{Protosec} we obtain an injection
$$i_A: A\times_T C' \hookrightarrow \JacS_{X'/C'}. $$
By Theorem  \ref{UnivFam} we obtain a morphism $\varphi':C' \to X'(d)$ 
which we can extend to $\varphi: C \to X(d)$, because $C$
is smooth. $\varphi$ is onto, since otherwise for any $t\in T$
the Jacobian $\JacS_{X'_t/C'_t}$ would be a product and $X_t \to C_t$
isotrivial. By Prop.\ \ref{genHELok} this contradicts $d \geq 3$.
The birational equivalence of $X$ and $C \times_{X(d)} S(A,d)$
is now nothing but the uniqueness assertion of Corollary
\ref{JtoC}.
\newline
To prove the uniqueness of $\varphi$, let $P\in C(\CC)$ be such 
that the fibre $F_P$ is smooth. The composition $m' \circ \varphi$ maps $P$
to the point in $A_3$ corresponding to the isomorphism
class of $\JacS_{F_P}$. The uniqueness now follows from
the injectivity of $m$ (Remark \ref{mrem}).
\par
{\it Step 3:} Consider again the curve 
$B'_n$ defined in the first step and let $h_n': S_n' \to B_n'$ be the pullback
of the universal family over $\mrC_3^{[n]}$
\newline
Via $\phi$ (see Remark \ref{mrem}) 
the factor group $G = \Gamma(d)/\Gamma(nd)$ acts without fixed points on 
$X'_n(d)$, on $M_3^{[n]}$ and thus on $B'_n$.
$G$ also acts on the universal family and hence on its
pullback to $B'_n$.
The quotient of the pullback by $G$, denoted by $h':S' \to B'$
does no longer depend on any choices (level structure, symplectic
basis) made above. We now check that the completion to a relatively 
minimal model $h(A,d): S(A,d) \to B(A,d)$ has the desired properties.
\newline
By construction $\JacS_{S_n'/B_n'}$ is globally isomorphic
to $J_n'\times_{X_n'(d)} B_n'$. This isomorphism is $G$-invariant and so
$\JacS_{S'/B'} \cong J' \times_{X'(d)} B'$.
This says that $h: X \to B$ has
maximal irregularity. Hence it will be a good prototype, once
we have shown that given a fibred family of surfaces $X \to C \to T$, 
the morphism $\varphi: C \to X(d)$ constructed in step $2$ factors 
via $B(A,d)$. 
\newline
After a suitable {\'e}tale base change $\tilde{C} \to C$ the
family $\tilde{X} \to \tilde{C}$ admits a level-$[n]$-structure
and the morphisms to $\tilde{C} \to \ol{\mrC_3}^{[n]}$ and $\tilde{C} \to
X_n(d)$ together with the smoothness of $\tilde{C}$ give a
morphism $\varphi_n: \tilde{C} \to B_n(A,d)$. The composition of
$\varphi_n$ with the projection $B_n(A,d) \to B(A,d)$ is independent
of the choice of the level structure and hence gives
the desired factorisation.
\newline
The uniqueness of $\varphi$ up to $\sigma$ follows from the
uniqueness statement for the family over $X(d)$.
\hfill $\Box$
\par
We now come back to the case $d=2$ and $\gamma=1$.
\par
{\bf Proof of Thm.~\ref{3protod2}:}
For each set of irregular cusps $S$ we can proceed as in section
\ref{secAbVar} to construct families of principally polarised
abelian varieties of dimension $3$
$$ j'(A,2)_S: J(A,d)_S \to X'(2)_S,$$
where of course $X'(2)_S \cong X'(2) \cong \PP^1 \sms\{0,1,\infty\}$
(as schemes, not as quotient stacks). We claim that over a connected
base $T$ each principally polarised abelian scheme $(J/T, \Theta_J)$
with an injection $A\times T \to J$ over $T$ satisfying $i_A^*(\Theta_J)
= \pr_A^*(\Theta_A)$
is the pullback of $j'(A,2)_S$ for a suitable $S$.
\newline
In fact we can construct locally on charts $T_i$
as in the proof of Thm.~\ref{UnivFam} sections $u^i_g$ and $u^i_{2g}$,
whose quotient in $V/V_A$ gives a morphism $T_i \to \HH$.
Similar sections  $u^k_g$ and $u^k_{2g}$ on another chart $T_k$
differ on $T_i \cap T_k$ form these by the action of $M \in \Gamma(2)$,
according to Lemma \ref{HomobasisW}. But the normalisation to
the upper half plane implies that $M$ belongs to one of the
complements of $\pm 1$.
\par
To prove the existence of the prototype we proceed as in 
Thm.~\ref{goodpro3}: we take
an auxiliary level structure and apply Corollary
\ref{JtoC} to obtain $h(A,2)_S$. 
\newline
For the prototype property we can also conclude as above,
once we have shown that we can distinguish the abelian
schemes $j(A,2)_S$ by its monodromy around the cusps.
Take w.l.o.g.\ the cusp at $\infty$. In the regular
and irregular case we have respectively
$$
\begin{array}{lcl}
{\rm Stab}_{\Gamma(2)_S}(\infty) = \langle 
\left(
\begin{array}{rr}
1 & 2 \\
0 & 1 \\
\end{array}
\right)
\rangle,
& \quad \text{hence} &
M(\gamma) = I_6 + E_{3,6} \\
{\rm Stab}_{\Gamma(2)_S}(\infty) = \langle
\left(
\begin{array}{rr}
-1 & 2 \\
0 & -1 \\
\end{array}
\right)
\rangle,
& \quad \text{hence} &
M(\gamma) =
\left(
\begin{array}{cccccc}
1 & 0 & 0 & 0 & 0 & 0 \\
0 & 1 & 0 & 0 & 0 & -1 \\
0 & 0 & -1 & 0 & 1 & -1 \\
0 & 0 & 0 & 1 & 0 & 0 \\
0 & 0 & 0 & 0 & 1 & 0 \\
0 & 0 & 0 & 0 & 0 & -1 \\
\end{array}
\right).
\end{array}
$$
These monodromy matrices are not conjugate in the symplectic group $Sp(3,\ZZ)$
and this proves the claim.
\hfill $\Box$
\par

\subsection*{Invariants of prototypes}
To calculate the invariants of these prototypes for $d \geq 3$ and for
irreducible $A$, we need to know, where $B(A,d) \to X(d)$
is actually ramified. 
\newline
We say that a point $P \in X(d)$ is {\em hyperelliptic}, 
if the image of the morphism $m: X(d) \to \ol{\mrC_3}$, which
extends $i^{-1} \circ m'$, is in the closure of the
hyperelliptic locus. This depends
of course on the abelian variety $A$. We use the same terminology
for the morphism $m_n: X_n(d) \to \ol{V_3^{[n]}}$.
We first give a proof of Proposition 3.13 in \cite{Na74}, 
which we have not been able to find in the literature.
\par
\begin{Prop} \label{monoprop}
Let $X \to B$ be a family of curves over a one-dimensional base $B$, 
smooth outside $P \in B$. If the monodromy around $P$ is unipotent, 
$X$ admits a stable model over $B$ (i.e.\ without base change).
\end{Prop}
\par
{\bf Proof:}
Let $B' = B \sms P$. By \cite{FC90} Theorem V.6.7 the
monodromy hypothesis implies that $\JacS_{X'/B'}$ has an
extension to a semi-abelian scheme $J \to B$, 
which is unique by \cite{FC90} Proposition I.2.7. Take
a cyclic cover $\tilde{B} \to B$ totally ramified over $P$, 
with covering group generated by $\sigma$, 
such that $\tilde{h}: \tilde{X} = X \times_B \tilde{B} \to
\tilde{B}$ admits a semistable
model. Thus $\JacS_{\tilde{X}/\tilde{B}} = J \times_B \tilde{B}$,
and $\sigma$ acts trivially on the fibre over $P$ of
this abelian scheme.
This means that $\sigma$ acts trivially on the fibre
of $\tilde{h}$ over $P$ and the claim follows.
\hfill $\Box$
\par
\begin{Lemma} \label{HEdurchnd}
$m_{n}$ does not factor in any complex neighbourhood of a hyperelliptic
point via $\ol{\mrC_3}^{[n]}$. The image of $B_n(A,d) \to \ol{\mrC_3}^{[n]}$
intersects the hyperelliptic locus transversally.
\end{Lemma}
\par
{\bf Proof:} Let $P \in B_n(A,d)$ lie over a hyperelliptic point of $X_n(d)$. 
Due to the injectivity of $m_n$ both claims follow, once we have
shown that the image of $\pr_2: B_n(A,d) \to \ol{\mrC_3}^{[n]}$ intersects
the hyperelliptic locus transversally at $\pr_2(P)$ and that
the tangent space of $\pr_2(B_n(A,d))$ at $\pr_2(P)$ is invariant
under the involution $\Sigma$.
\newline
Take a punctured neighbourhood $U'$ of $P$ and let $h_{U'}$ be
the pullback of $h(A,d)$ to $U'$. By construction the
Jacobian of $h_{U'}$ is $J_{U'} = J(A,d) \times_{X(d)} U'$ and
the monodromy around $P$ is unipotent by Lemma \ref{MonodrLemma}.
By the above proposition we can extend $h_{U'}$ to a stable family of
curves over $U = U' \cup \{P\}$. Let $C$ be the fibre
of $h_U$ over $P$. We distinguish whether $P$
maps to $X'_n(d)$ (first case) or to a cusp (second case).
\newline
Denote by $p: J_{U'} \to A^\vee \times U'$ the dual
of the inclusion $i_A$ of the fixed part. Consider
how the curve $p(C)$ is deformed on $A^\vee$ instead
of just considering abstract deformations.
In the first case the geometric genus of $p(C)$ is $3$,
because otherwise $d=2$. As the canonical bundle on $A^\vee$
is trivial, first order deformations of the normalisation
of $p(C)$ (or equivalently: of the pair $(C,p)$) are
parametrised by $H^0(C, K_C)$ (see \cite{Ta84}, 
section before Lemma 1.5). 
The Kodaira-Spencer-mapping
hence induces the following commuting triangle, equivariant
with respect to the action of the hyperelliptic involution:
$$
\xymatrix{
T_{B_n(A,d),P} \ar[rr]^{\kappa_{h,P}} \ar[dr]^{\kappa_{P}} & & H^1(C,T_C) \\
& H^0(C,K_C)  \ar[ur]^\delta \\
}$$
$\delta$ stems from the connecting homomorphism of
$$ 0 \to T_C \to T_A|_C \to K_C \to 0.$$
The map $\kappa_P$ is not zero because $J_U$
is not a trivial deformation of $\JacS_C$ in any 
neighbourhood of $P$, as one can see already by using
the family $J(z) \to \HH$.  The hyperelliptic
involution acts on $H^0(C,K_C)$ as $(-1)$, and so
the image of $\kappa_{h,P}$ lies in the ($1$-dimensional)
eigenspace of $(-1)$ of $H^1(C,T_C)$.
\newline
The other case works essentially the same way: If $\pi:\tilde{C}
\to C$ is the normalisation, first order deformations of
the pair $(\tilde{C}, p \circ \pi)$ are parametrised by
$H^0(\tilde{C}, K_{\tilde{C}})$ (see \cite{Ta84} Lemma 1.5
and Remark 1.6). We thereby use the fact that $C$ is stable
and therefore the ramification divisor of $\pi$ is trivial.
On $H^0(\tilde{C}, K_{\tilde{C}})$ the hyperelliptic involution
still acts as $(-1)$ and we conclude as above. 
\hfill $\Box$
\par 
\begin{Cor} \label{covramhyp}
The twofold covering $B(A,d) \to X(d)$ is ramified 
precisely over the hyperelliptic points of $X(d)$.
\end{Cor}
\par
Let $t(d)$ be the number of cusps of $X(d)$ and
$s_g(A,d)$ the number of points where $J(A,d)$ is proper, but reducible for
$g=2$ or $g=3$ respectively. Defining 
$$ \Delta_d = \frac{d^2}{24}\prod_{p\mid d} (1-\frac{1}{p^2}),$$
we know from \cite{Sh71}:
$$ g(X(d)) = (d-6)\Delta_d + 1$$
$$ t(d) = 12\Delta_d.$$
\par
We first sum up Xiao's results (\cite{Xi85}) for $g=2$:
\begin{Cor} \label{Xiaoresults}
For $g=2$, $s_2(d) = s_2(A,d)$ does not depend on $A$. More
precisely:
$$\begin{array}{lllll}
s(d) &=& (5d-6)\Delta_d \\
c_2(S(A,d)) &=& s(d) + t(d) + 4g(X(d)) - 4 &=& 
(9d-18)\Delta_d\\
\chi({\cal O}_{S(A,d)}) &=& 2g(X(d)) - 2 + \frac{1}{2}t(d) 
&=& (2d-6) \Delta_d\\
K_{S(A,d)}^2 &= & 6\chi({\cal O}_{S(A,d)}) + 3g(X(d)) - 3
&=& (15d-54)\Delta_d\\
\end{array} $$ 
The singular fibres of surfaces in $\mfS_{2,b}^{\mi}$
are two elliptic curves attached to each other at a node, if the base point 
maps in $X(d)$ to one of the $s_2(d)$ points corresponding to a proper
but reducible Jacobian or an elliptic curve with a node, if the base
point maps to a cusp of $X(d)$.
\end{Cor}
\par
{\bf Sketch of Proof:}
First, in order to prove Theorem \ref{goodpro2}, the family
of polarisation divisors of $j'(A,d)$ is already the family
of curves $h'(A,d)$. Hence we have $\JacS_{S(A,d)/B(A,d)} \cong J(A,d)$
globally and no trouble with Torelli. The proof of the prototype
property and the uniqueness of $\varphi$ works as above.
\newline
The irreducible fibres can be recognised by calculating their
monodromy as in Lemma \ref{MonodrLemma} and comparing with the list in
\cite{NU73}. To calculate the invariants, the techniques are
similar to the ones below.
\hfill $\Box$
\par
In the case $g=3$ we obtain:
\begin{Cor} \label{fibreresults}
For $d\geq 3$ and $A$ irreducible the genus of $B(A,d)$ and 
$s_3(d) = s_3(A,d)$ does not depend on $A$. Moreover we have: 
$$\begin{array}{lllll}
s(d) &=& 0 \\
g(B(A,d)) &=& (20d-36)\Delta_d + 1 \\
c_2(S(A,d)) &=& (160d-264)\Delta_d\\
\chi({\cal O}_{S(A,d)}) &=& (42d-72) \Delta_d\\
K_{S(A,d)}^2 &=& (344d-600)\Delta_d\\
\end{array} $$ 
The only singular fibres of surfaces in $\mfS_{3,b}^{\mi}$
are curves of genus $2$ with one node or a smooth curve
of genus $2$ together with a smooth $\PP^1$, which have
two transversal intersections.
\newline
In particular $S(A,d) \to B(A,d)$ is semistable. The index of
$S(A,d)$ is $\tau = (8d-24)\Delta_d$, hence positive for $d>3$.
\end{Cor}
\par
{\bf Proof:} First, $s_3(A,d) = 0$ was shown in Lemma \ref{3noredlemma}.
For simplicity, we now drop $(A,d)$ from the notation.
\newline
Secondly, we examine the fibres of $h: S \to B$. 
By Cor.\ \ref{covramhyp} the map $B_n \to X_n(d)$ is of degree two
with the hyperelliptic points as ramification locus.
By Lemma \ref{MonodrLemma}
the monodromy around the cusps of $B(A,d)$ is $I_6 + 2E_{3,6}$
or $I_6 + E_{3,6}$, depending on whether the cusp is hyperelliptic
or not. Thus by Proposition \ref{monoprop} 
$h$ is semi-stable and the non-smooth fibres
lie over the cusps. The period matrix calculated 
before  Lemma \ref{MonodrLemma}
implies that the Jacobian of these fibres  is
an extension of an abelian surface by a torus. To prove the
assertion, we must exclude that these fibres consist of 
a genus $2$ curve attached to a $\PP^1$ with a self-intersection.
But in this case, the period matrix would have a block structure.
\newline
Thirdly we calculate the number of hyperelliptic fibres of $h$. 
For this purpose we use the following equation (\cite{HM98} (3.165)):
$$ H = 18\lambda -2 \delta_0 -3 \delta_1.$$
$H$ is number of hyperelliptic fibres, $\lambda$ is the
degree $h_*\omega_{S/B}$ and $\delta_i$ is the number of
curves belonging to the boundary components $\Delta_i$
(both $H$ and $\delta_i$ have to be counted with multiplicity).
Calling $j_{B}: J_{B} \to B$ the pullback of $j(A,d)$ to $B$, we have
$$\deg {h}_* \omega_{S/B} = 
-\deg (\bigwedge^3 R^1 {j_{B}}_* {\cal O}_{J_{B}})
= -2 \deg (\bigwedge^3 R^1 j(A,d)_* {\cal O}_{J(A,d)}).$$
We can evaluate the right-hand side as in \cite{Xi85} Thm.~3.10
using modular forms and we obtain
$$\deg (\bigwedge^3 R^1 j(A,d)_* {\cal O}_{J(A,d)}) = -(g(X(d))- 1 + t(d)/2) =
 -d\Delta_d.$$
The monodromy implies $\delta_0 = 2t(d)$ and we have seen
that $\delta_1 = s_3(d) = 0$. We hence conclude
$$H = (36d-48)\Delta_d$$
and by Lemma \ref{HEdurchnd} this is the number of
hyperelliptic fibres of $h$. This enables us to compute
the genus of $B=B(A,d)$.
\newline
To determine $c_2(S(A,d))$ we let $F$ be a generic fibre of $h$ and 
use the formula
$$\sum_{F' \text{singular}} (\chi_{\rm top}(F') - \chi_{\rm top}(F))
 = c_2(S(A,d)) - \chi_{\rm top}(F) \cdot \chi_{\rm top}(B(A,d)).$$ 
Note that the non-hyperelliptic singular fibres contribute by one
to $\chi_{\rm top}(F)$, the hyperelliptic fibres by two.
\newline 
Finally, to determine $\chi({\cal O}_S)$ remember that (by the
Leray spectral sequence and Riemann-Roch)
$$ \chi({\cal O}_S) = \chi({\cal O}_{B}) 
- \chi(R^1 h_*  {\cal O}_{S}) = 2(g(B) -1) - \deg(\bigwedge^3 R^1 h_*  
{\cal O}_{S}).$$
Using the above calculation, this degree is $-2d\Delta_d$ and
the claim follows.
\hfill $\Box$
\par
\begin{Rem}{\rm If $\tilde{X} \to \tilde{B}$ is the
minimal model of a surface in $\mfS_{3,\tilde{b}}^{\mi}(\CC)$
with a prototype base change $\varphi: \tilde{B} \to B(A,d)$
of degree $n$, we have 
$$ K_{\tilde{X}}^2 = nK_{S(A,d)}^2 + 8(\tilde{b} -1 + nb-n).$$
If $P \in \tilde{B}$ lies over a cusp of $B(A,d)$ and ramifies
to the order $e$, topology implies that
$$ \chi_{\rm top}(F_P) - \chi_{\rm top}(F) = e $$
and by the above formulae we can determine completely the
possible invariants of surfaces in $\mfS_{3,\tilde{b}}^{\mi}(\CC)$.
}\end{Rem}
\par
\begin{Rem}{\rm Let $S$ be a projective  surface of general type.
Assume that $S$ is of maximal  Albanese
dimension and   the  canonical map $\Phi_S$
is composite with a pencil. Note that then $q(S)=2$.
Suppose the generic fibre of the fibration $f: S\to C$   
associated with $\Phi_S$ has genus $3$. \newline
By Thm.~\ref{goodpro3}  and Cor.~\ref{fibreresults}
the associated degree of $f$ is $2$.
Hence such surfaces are classified by Theorem 2.8.
Explicit examples of such surfaces  are given in \cite{Be79} Ex.\ 2
and \cite{Xi87b} Ex.\ 3. 
}\end{Rem}

\section{The components of the moduli space for maximally
irregularly fibred surfaces} \label{secMSpace}

\subsubsection*{Families of fibred surfaces and fibred families of surfaces}
\label{ffam}

Theorems by Beauville and Siu (for the base curve, \cite{Be91},
\cite{Siu87})) and Catanese (for the fibre, \cite{Ca00}) state that
the property of having a fibration of type $(g,b)$ is deformation
invariant, if both $b \geq 2$ and $g \geq 2$. We call this
the fibre genus condition $(FGC)$ and base genus condition $(BGC)$
respectively, and {\em we restrict ourselves from now 
on to such surfaces.} Due to these theorems we can state:
\par
\begin{Defi} We denote by $\mfS_{g,b}(\cdot)$ the (open and closed)
subfunctor of
$\mfS(\cdot)$ parametrising families of surfaces of general type
such that each fibre admits a fibration of type $(g,b)$. We denote
by $N_{g,b}$ the union of the corresponding components of the moduli
space.
\end{Defi}
\par
For the analysis of the moduli space of surfaces, the 
question remains, 
whether a family $X \to T$ of fibred surfaces in $\mfS_{g,b}(T)$
is, what we call a {\em fibred family of surfaces}, that is if
there is a family of curves $B \to T$ such that $X \to B \to T$ induces
the fibrations of type $(g,b)$. For regular fibrations, the
Albanese mapping gives an affirmative answer.
We now give a sufficient condition for families of irregularly 
fibred surfaces to be fibred
in families and a numerical criterion for testing the condition.
\par
\begin{Defi}
A family of surfaces $X \to T$ is said to have {\em at most one
fibration}, if given two fibrations $h_i: X \to B_i$ over $T$,
there is a $T$-isomorphism $\psi: B_1 \to B_2$ satisfying 
$\psi \circ h_1 = h_2$.
\newline
The subfunctor $\mfS_{g,b}(\cdot)$ of $\mfS(\cdot)$ 
is said to have this property, 
if it is fulfilled for each $T$ and each surface in $\mfS_{g,b}(T)$.  
\end{Defi} 
\par
\par
\begin{Thm} \label{FibinFam}
If $\mfS_{g,b}(\cdot)$ 
parametrises irregularly fibred surfaces of type $(g,b)$ with
at most one fibration, then there is a natural transformation
$\mfS_{g,b} \to \mfC_b$ which (on complex points)  
sends each surface to the base of the fibration.
Hence we have a natural morphism
$$ N_{g,b} \to M_b$$
between the corresponding coarse moduli spaces.
\end{Thm}
\par
{\bf Proof:} It is well known (\cite{Ca91} Theorem 4.9 or
\cite{Ser92} Claim 5.1) that, given a fibration $h:X \to B$,
one has a surjection 
of deformation functors $\Def_h \to \Def_X$,
provided that the cohomology group $H^0(R^1 h_* {\cal O}_X \otimes
T_X)$ vanishes. This vanishing property follows from relative duality and
a theorem of Fujita (see loc.\ cit.). This solves
the problem for a local Artinian base. If the base
is a complete local ring $R$, Grothendieck's
algebraisation theorem ([EGA] III th{\'e}or{\`e}me 5.4.5) 
implies that the deformations over the Artinian quotients
stem from  a deformation of $B$ over $R$. 
\newline
For the general case, we cover the base $T$ by the spectra of
the completions of its local rings. As we suppose the base to be Noetherian, 
finitely many will be sufficient and we call the union of
these schemes $\tilde{T}$. What we need is an (fppf-) descent
datum on the base curve $\tilde{B} \to \tilde{T}$ obtained
by the above deformation techniques. By $(BGC)$ $\tilde{B}$
comes along with an ample canonical sheaf and the descent
datum will automatically be effective (\cite{BLR90} Theorem 6.1.7).
\newline
$\tilde{X} = X \times_T \tilde{T}$ has 
a natural morphism $\varphi: \pr_1^* \tilde{X} \to 
\pr_2^* \tilde{X} $, 
where $\pr_i$ are the usual projections $ \tilde{T} \times_T \tilde{T} 
\to \tilde{T}$. Given the fibration $\tilde{h}: \tilde{X} \to
\tilde{B}$, we can apply the hypothesis to
$\pr_1^*(\tilde{h})$ and $\pr_1^*(\tilde{h}) \circ \psi$ 
to obtain $\psi: \pr_1^* \tilde{B} \to 
\pr_2^* \tilde{B} $. 
The map $\psi $ obviously satisfies the cocycle condition, 
because $\varphi$ does and because $\tilde{h}$ is surjective.
\hfill $\Box$
\par
\begin{Crit} \label{onefibcrit}
A family of surfaces $X \to T \in \mfS_{g,b}(T)$, whose
fibres $X_t$ satisfy 
$$ K^2_{X_t} > 4 (g-1)^2,$$
has at most one fibration of type $(g,b)$.
\end{Crit}
\par
{\bf Proof:} In case $T = \Spc \CC$ this is \cite{Xi85} Proposition 6.4.
In general, suppose there are two fibrations $X \to B_i \to T$ 
($i=1,2$). Consider
the product morphism $h: X \to B_1 \times_T B_2$. The image $Y$ of $h$ is
flat over $T$ by the criterion in [EGA] IV, 11.3.11. Hence the image is
a family of curves. Consider the projections $\pr_i: Y \to B_i$. 
They are finite and isomorphisms for each closed point $t \in T$.
Hence $\pr_i$ is {\'e}tale (this is an open condition!) and radicial
(by the case $T= \Spc \CC$ and [EGA] I, 3.7.1). Now [EGA] IV, 17.9.1
proves that $\pr_i$ is an isomorphism and the criterion follows.
\hfill $\Box$
\par
The families of fibred surfaces which we investigate here satisfy
the hypothesis of
this criterion:  Let $\mfS_{g,b}^{\mi}(\cdot)$ be the subfunctor
of $\mfS_{g,b}$ parametrising fibred surfaces with maximally irregularity.
Accordingly, $N_{g,b}^{\mi}$ denotes the union of the corresponding components
of the moduli space. 
\par
\begin{Lemma} \label{FibfamAppl}
A family of fibred surfaces 
in $\mfS_{g,b}^{\mi}(T)$ with $g=2$ or $g=3$,
which is not isotrivial,
has at most one fibration and is hence a fibred family.
\end{Lemma}
\par
{\bf Proof:} To apply Criterion \ref{onefibcrit} to 
$X \in \mfS_{g,b}^{\mi}(\CC)$, 
the Arakelov inequality $$K_X^2 \geq 8(b-1)(g-1)$$ (see \cite{Be82})
is sufficient for $g=2$ thanks to $(BGC)$. For $g=3$ we have to 
exclude the case $b=2$ and $K_X^2 = 16$. Using the well known relation
of $K_X^2, \chi({\cal O}_X)$ and the slope $\lambda$ (see \cite{Xi87})
$$ K_X^2 = \lambda \echar_X) + (8-\lambda)(b-1)(g-1),$$
we deduce $\chi({\cal O}_X) = 2$ and this implies that $X$ is
a fibre bundle (see again \cite{Be82}), hence isotrivial.
\hfill $\Box$
\par
\subsubsection*{Theorem \ref{UnivFam} in families} 
We need a relative version of this theorem obtained
by letting $Z$ vary in $\HH_{g-1}$. 
Dividing out the $\Gamma(d)$-action on the family 
$J(z,Z) \to \HH \times \HH_{g-1}$ constructed in Section $2$
(there with fixed $Z$), we obtain a family of principally 
polarised abelian varieties
over $X'(d) \times \HH_{g-1}$. For technical reasons
we do not take the quotient of $\HH_{g-1}$ by the whole symplectic
group but fix a subgroup $G$, such that
$\HH_{g-1}/G$ is a moduli space for abelian varieties
with a polarisation of type $\delta$ and a level-$[n]$-structure. 
The quotient by $G$,
$$ j'(d): J'(d) \to X'(d) \times \mrAV_{g-1,\delta}^{[n]},$$
is a family of polarised abelian varieties over the product of the
modular curve and the moduli space of abelian varieties with
level-$[n]$-structure. The same proof of Thm.\ \ref{UnivFam} yields: 
\par
\begin{Cor} \label{Unifam}
$j'(d)$ is the universal family of $g$-dimensional, principally polarised
abelian schemes with an injection of an $(g\!-\!1)$-dimensional 
abelian scheme $A$ endowed with level-$[n]$-structure, 
provided that $d \geq 3$.
\end{Cor}

\subsubsection*{Components of the Gieseker moduli space}
The results in the previous section suggest to describe
$N_{g,b}^{\mi}$ by the fixed part of the fibration
plus the morphisms of the base curve to $B(A,d)$. For that
purpose, let $\mfAV_{g-1,\delta}(\cdot)^f$ denote the subfunctor
parametrising abelian schemes occurring as fixed
parts of maximally irregular fibrations of type $(g,b)$.
According to the Theorems \ref{goodpro2} and \ref{goodpro3}
this means no restriction for $g=2$
and no principally polarised $1$-dimensional abelian subvarieties
for $g=3$. We denote by $\mrAV_{g-1,\delta}^f$ the
corresponding subscheme of $\mrAV_{g-1,\delta}$.
\par
For a fixed curve $D$ we need the following functor for
coverings:
$$ \mfC_{b}(D,m)(T) = \{ (B/T, \varphi: B \to D) \, 
{\rm where} \,B \in \mfC_{b}(T), 
\, \deg(\varphi) = m \}. $$
This is an open subfunctor of a stable mappings functor,
that (according to \cite{FuPe95}) has a coarse moduli space
$M_b(D,m) \subseteq \ol{M_b}(D,m)$. 
\par
We can now use the prototype to prove a structure theorem for
the moduli spaces $N_{g,b}^{m.i.}$:
\begin{Thm} \label{strucresult}
The moduli space for maximally irregularly fibred surfaces
of type $(g,b)$ ($g \in \{2,3\}$) 
with associated degree $d \geq 3$ decomposes
into components according to $d$ and the degree $m$ of
the morphism $\varphi: B \to B(A,d)$ induced by the fibration $X \to B$.
\newline
In case $g=2$ each such component is isomorphic to
$$ M_b(X(d), m) \times \mrAV_1$$
and has dimension $2b-2-m(2g(X(d)) -2) + 1$.
\newline
In case $g=3$ each such component admits a morphism
to $\mrAV_{2,\delta}^f$ with the fibre over $A$ isomorphic
to $M_b(B(A,d),m)/\sigma$, where $\sigma$ comes from the
involution of $B(A,d) \to X(d)$. The dimension of the
moduli space is $2b-m(2g(B(A,d)) -2) + 3$. 
\end{Thm}
\par
{\bf Proof:} Recall that by Lemma \ref{FibfamAppl}
each family of surfaces $X \to T$ in
question is a fibred family $X \to B \to T$. The fixed part of this
fibration induces a morphism $\psi: T \to  \mrAV_{g-1,\delta}^f$.
After an {\'e}tale base change $\tilde{T} \to T$ we may suppose that
the fixed part has a level-$[n]$-structure inducing 
$\tilde{\psi}: \tilde{T} \to  (\mrAV_{g-1,\delta}^{[n]})^f$.
We use the level-structure for a moment to apply Corollary \ref{Unifam}.
The morphism $i_A: A_{B'} = A \times_T B' \to \JacS_{X'/B'}$ 
constructed at the beginning of Section \ref{Protosec} gives 
$\varphi': B' \times_T \tilde{T} \to X'(d)$. As $B$
is smooth $\varphi'$ extends to $B \times_T \tilde{T}$. By the
universal property of $j'(d)$ this morphism is unique and therefore
glues with the descent data on $\tilde{T}$ to give $\varphi:B \to X(d)$.
In view of Theorem \ref{goodpro3} we have shown, that the moduli 
functor admits a natural transformation to the product 
$\mfC_{b}(X(d),m) \times \mfAV_{1}$ in case $g=2$
and in case $g=3$ that it has a natural transformation to
$\mfAV_{2,\delta}^f$, whose fibres over $A$ are in 
$\mfC_{b}(B(A,d),m)$.
\newline
Conversely, given $(A,\Theta_A) \in \mrAV_{g-1,\delta}^f (\CC)$ and 
$\varphi \in M_b(B(A,d),m)$, the property 'good prototype' gives
a fibred surface with the prescribed $A$, $d$, $m$ and $\varphi$. If
two fibrations $h_i: X_i \to B_i$ associated to $(\phi_i, A_i)$
are isomorphic, the fixed parts $A_i$ are necessarily isomorphic as
polarised abelian schemes.
As the fibrations are unique by Lemma \ref{FibfamAppl},
we obtain an isomorphism
$\iota: B_1 \to B_2$ and it remains to show that
$\phi_2 \circ \iota = \phi_1$, resp.\ up to $\sigma$ in
the case $g=3$. If this was wrong at $P \in X'(d)$, compare
the Jacobians over the preimage in $B_i$ to obtain 
a contradiction
to the injectivity of $m'$ (see Rem.~\ref{mrem}).
\newline
The maximality among the coarsely representing spaces
in the case $g=2$ is easily checked by redoing this argument
with a family $A/T$ instead of $A$ . 
In case $g=3$, we only claimed such a statement for fixed
$A \in \mrAV^f_{2,\delta}$ and this was already done in the
proof of theorem \ref{goodpro3}.
\hfill $\Box$
\par
The proof shows that for any $g \geq 2$ and any $d \geq 2$ 
we have a natural transformation
$$ \eta: \mfS_{g,b} \to \mfC_{b}(X(d),m) 
 \times \mfAV^f_{g-1,\delta}$$ 
and therefore a morphism
$ \eta: N_{g,b} \to \mrC_{b}(X(d),m) \times \mrAV^f_{g-1,\delta}$,
which is injective for $d \geq 3$ and at most $4:1$ 
for $d = 2$ and $\gamma =1$, provided that all families
are fibred families (e.g.\ for $b$ big enough). 
\par
\begin{Rem} {\rm For $g \geq 4$ the crux of the matters is to determine what 
$\mrAV^f_{g-1,\delta}$ in fact is, i.e.\ to determine the abelian 
varieties of dimension $g-1$, that occur as fixed parts of 
one-dimensional families of Jacobians of dimension $g$. 
\newline
Pirola showed in \cite{Pi92} that $\mrAV^f_{g-1,\delta}$ is non-empty
for $g=4$ but it is an interesting question to decide
whether it is empty or not for $5 \leq g \leq 7$. }
\end{Rem}
\par

\par
Martin M{\"o}ller: Universit{\"a}t Essen, FB 6 (Mathematik) \newline 
45117 Essen, Germany \newline
e-mail: martin.moeller@uni-essen.de \newline

\end{document}